\documentclass{article}

\usepackage{verbatim}
\usepackage{graphicx}
\usepackage{amssymb}

\newcommand{\Beq}{\begin{equation}}
\newcommand{\Eeq}{\end{equation}}

\begin{document} 

\title{Mixed connectivity of Cartesian graph products and bundles%
       \thanks{This work was supported in part by the Slovenian research agency.} 
}
 
\author{ Rija Erve\v{s},\\ 
FCE, University of Maribor,\\ 
Smetanova 17, Maribor 2000, Slovenia\\
{\tt rija.erves@uni-mb.si}  
\and Janez \v{Z}erovnik\\
FME, University of Ljubljana,\\ 
A\v sker\v ceva 6, \\
SI-1000 Ljubljana, Slovenia\\
and \\
Institute of Mathematics, Physics and Mechanics,\\
Ljubljana, Slovenia\\
{\tt janez.zerovnik@fs.uni-lj.si, janez.zerovnik@imfm.si }
}
\date{\today}

\maketitle 

\begin{abstract}
Mixed connectivity is a generalization of  vertex and edge connectivity. 
 A graph is $(p,0)$-connected, $p>0$, if the graph remains connected after removal of any $p-1$ vertices. 
A graph   is $(p,q)$-connected, $p\geq 0$, $q>0$, 
if it remains connected after removal of any $p$ vertices and any $q-1$ edges. 
 Cartesian graph bundles are graphs that  generalize both covering graphs  and Cartesian graph products. 
It is shown that if  graph $F$ is $(p_{F},q_{F})$-connected and graph $B$ is $(p_{B},q_{B})$-connected, 
then Cartesian graph bundle $G$ with fibre $F$ over the base graph $B$ is $(p_{F}+p_{B},q_{F}+q_{B})$-connected.
Furthermore, if $q_{F},q_{B}>0$, then  $G$ is also $(p_{F}+p_{B}+1,q_{F}+q_{B}-1)$-connected.
Finally, let graphs $G_i, i=1,\dots,n,$ be $(p_i,q_i)$-connected and let $k$ be the number of graphs with $q_i>0$. 
The Cartesian graph product $G=G_1\Box G_2\Box \dots \Box G_n$ is
$(\sum p_i,\sum q_i)$-connected, and,  for $ k\geq 1$, it is also  $(\sum p_i+k-1,\sum q_i-k+1)$-connected. 
\end{abstract}
 

\newtheorem{theorem}{Theorem}[section]
\newtheorem{corollary}[theorem]{Corollary}
\newtheorem{proposition}[theorem]{Proposition}
\newtheorem{lemma}[theorem]{Lemma}
\newtheorem{remark}[theorem]{Remark}
\newtheorem{definition}[theorem]{Definition}
\newtheorem{question}[theorem]{Question}
\newtheorem{example}[theorem]{Example}
\newtheorem{conjecture}[theorem]{Conjecture}

\newcommand{\proof}{\noindent{\it Proof:~}} 
\newcommand{\qed}{{\hfill$\Box$}} 

\section{Introduction}

 Graph products  and bundles are among frequently  studied interconnection network topologies.
For example the meshes, tori, hypercubes and some of their generalizations are Cartesian products.
It is less known that  some well-known topologies 
are  Cartesian graph bundles, i.e. some twisted hypercubes \cite{cull,efe} and 
multiplicative circulant graphs \cite{ivan}.  
Other graph products, sometimes under different names, have been studied 
as interesting  network topologies \cite{bermond,Xavier,ivan}.

In the design of large interconnection networks several factors have to be taken into account.
A usual constraint is that each  processor can be connected to a limited number of other processors and 
the delays in communication must not be too long. 
Furthermore, an interconnection network should be fault tolerant, because practical communication networks are exposed 
to failures of network components. 
Both failures of nodes and failures of connections between them happen and it is desirable that a network is robust in 
the sense that a limited number of failures does not break down the whole system.
A lot of work has been done on various aspects of network fault tolerance,
see for example the survey \cite{Bermond} and more recent papers \cite{hung,sun,yin}.
In particular the fault diameter with faulty vertices which was first studied in \cite{1} and the edge fault diameter has  
been determined for many important networks recently \cite{day,4,7,0}. 
In particular, the (vertex) fault diameter and the edge fault diameter of Cartesian graph products and Cartesian graph bundles 
was studied recently \cite{zer-ban,zbedge,ZB,erves}.
Usually either only edge faults or only vertex faults are considered, 
while the case when both edges and vertices may be faulty is studied rarely. 
For example, \cite{hung,sun} consider Hamiltonian properties assuming a combination of vertex and edge faults. 
In recent work on  fault diameter of Cartesian graph products and bundles \cite{zer-ban,zbedge,ZB,erves}, 
analogous results were found for both fault diameter and edge fault diameter. However, the proofs
for vertex and edge faults in \cite{zer-ban,zbedge,ZB,erves} are independent, 
and our effort to see how results in one case may imply the others was  not successful.
A natural question is whether it is possible to design  a uniform theory that would enable unified proofs or 
provide tools to translate results for one type of faults to the other. 
It is therefore of interest to study general relationships between invariants under vertex and edge faults.
Some basic results on edge, vertex and mixed fault diameters for general graphs appear in \cite{3fd}.
In order to study the fault diameters of graph products and bundles  under mixed faults, it is important to 
understand the generalized connectivities.

Here we study  mixed connectivity which generalizes both vertex and edge connectivity.
It is known  that Cartesian graph bundle  with fibre $F$ over the base graph $B$ is  
$(\kappa(F)+\kappa(B))$-connected and $(\lambda(F)+\lambda(B))$-edge connected \cite{zer-ban,erves}.
In this paper we generalize these results to mixed connectivity of Cartesian graph bundles.
More precisely,  
assuming that the fibre $F$ is $(p_{F},q_{F})$-connected and the base graph  $B$ is $(p_{B},q_{B})$-connected, 
then Cartesian graph bundle $G$ with fibre $F$ over the base graph $B$ is
$(p_{F}+p_{B},q_{F}+q_{B})$-connected.
Furthermore, if $q_{F},q_{B}>0$, then the Cartesian graph bundle is also $(p_{F}+p_{B}+1,q_{F}+q_{B}-1)$-connected.
As a corollary,  mixed connectivity of the Cartesian product of finite number of factors is given.

The rest of the paper is organized as follows.
In the next section  some general definitions are given, 
and in Section 3  the mixed connectivity  is defined and some basic facts are observed.
Section 4 recalls definition of graph bundles and states the main result (Theorem \ref{mpsv}) which is proved in the last section.



\section{Preliminaries}
 
 Here we only recall some basic definitions to fix the notation, 
for other standard  notions not defined here we adopt the usual terminology (see for example \cite{m}).
 A {\em simple graph}  $G=(V,E)$  is determined by a {\em vertex set} $V=V(G)$ 
and a set  $E=E(G)$ of (unordered) pairs of vertices, called the set  of {\em edges}.
As usual, we  will use the short notation $uv$ for edge $\{u,v\}$.
For an edge $e=uv$ we call $u$ and $v$ its {\em endpoints}.
It is convenient to consider the union of  \emph{elements} of a graph, $S(G)= V(G) \cup E(G)$.
Given $X \subseteq S(G)$ then $S(G) \setminus X$ is a subset of elements of $G$.
However, note that in general  $S(G) \setminus X$ may not induce a graph.
As we need notation for subgraphs with some missing (faulty) elements, 
we will formally  define $G\setminus X$, the subgraph of $G$ after deletion of $X$, as follows:

\begin{definition}
Let $X \subseteq S(G)$, and $X=X_E \cup X_V$, where $X_E \subseteq E(G)$ and $X_V \subseteq V(G)$.
Then $G\setminus X$ is the subgraph of $(V(G), E(G)\setminus X_E)$ induced on vertex set $V(G)\setminus X_V$. 
\end{definition}
A {\em walk} between $x$ and $y$ is a sequence of vertices  and edges
$v_0,$ $e_1,$ $v_1,$ $e_2,$ $v_2,$ $\dots,$ $v_{k-1},$ $e_k,$ $v_k$
where $x=v_0$, $y=v_k$, and  $e_i = v_{i-1}v_i$ for each $i$.
The {\em length} of a walk $W$, denoted by $\ell (W)$, is the number of edges in $W$. 
A walk with all vertices distinct is called a {\em path},
and the vertices $v_0$ and $v_k$ are called the {\em endpoints} of the  path.
A path $P$ in $G$, defined by a sequence 
$x=v_0,e_1,v_1,e_2,v_2,\dots,v_{k-1},e_k,v_k=y$ can alternatively be seen as a subgraph of $G$
with $V(P) =\{v_0,v_1,v_2,\dots,v_k\}$ and $E(P) =\{e_1,e_2,\dots,e_k\}$.
Note that the reverse sequence gives rise to the same subgraph. 
Hence we use $P$ for a path either from $x$ to $y$ or from $y$ to $x$.
A graph is {\em connected} if there is a path between each pair of vertices, and is {\em disconnected} otherwise.

The  {\em connectivity} (or {\em vertex connectivity}) of a connected graph $G$, $\kappa (G)$, is the minimum cardinality over all 
vertex-separating sets in $G$. 
As the complete graph $K_n$ has no vertex-separating sets, we define $\kappa (K_n)=n-1$. 
We say that $G$ is  {\em $k$-connected} (or {\em $k$-vertex connected}) for any $k \leq \kappa (G)$.
The  {\em edge connectivity} of a connected graph $G$, $\lambda (G)$, is the minimum cardinality over all edge-separating sets in $G$.
A graph $G$ is said to be {\em $k$-edge connected}  for any $k \leq \lambda (G)$.
In other words, the edge connectivity $\lambda (G)$ of a connected graph $G$ is the smallest number of edges whose removal disconnects $G$, 
and the (vertex) connectivity $\kappa (G)$ of a connected graph $G$ (other than a complete graph) is the smallest number of vertices whose 
removal disconnects $G$. 
It is well-known that
(see, for example, \cite{m}, page 224)
$\kappa (G) \leq \lambda (G) \leq \delta_G, $
where $\delta_G$ is the smallest vertex degree of $G$. 
Thus if a graph $G$ is $k$-connected, then it is also $k$-edge connected. The reverse does not hold in general.
For later reference recall  that by  Menger's theorems (see, for example, \cite{m}, pages 230,234) we know that 
in a $k$-connected graph $G$  there are at least $k$ vertex disjoint paths between any two vertices in $G$, 
and if   $G$ is $k$-edge connected then there are at least $k$ edge disjoint paths between any two vertices in $G$.

\section{Mixed connectivity}

Considering a graph with faulty vertices and faulty edges at the same time we can generalize both vertex and edge connectivity.
We start with an observation that can be proved easily.

\begin{proposition}
\label{T1} Let $p,q>0$. 
If a graph $G$ remains connected after removal of any $p$ vertices and any $q-1$ edges, 
then $G$ also remains connected after removal of any $p-1$ vertices and any $q$ edges.
\end{proposition}

\proof{Proof.} Let $H$ be any subgraph of $G$ after removal of $p-1$ vertices and $q-1$ edges. 
By assumption, the graph $H$ remains connected after removal of any vertex.
Hence $H$ is $2$-connected 
which implies that $H$ is $2$-edge connected.
In other words, $H$ remains connected after removal of any edge 
which in turn implies that $G$ remains connected after removal of any $p-1$ vertices and any $q$ edges.
\qed

If a graph $G$ remains connected after removal of any $p$ vertices (and any $q-1$ edges), then 
$p<\kappa(G)$. By repeated application of Proposition \ref{T1}, the graph $G$ also remains 
connected after removal of any $p+q-1$ edges, hence $p+q\leq\lambda(G)$. 
Now we  formally define mixed-connectivity.

{\begin{definition}
%
\hfill\break\noindent (1)  
Let $p>0$. 
Graph $G$ is \emph{$(p,0)$-connected}, if $G$ remains connected after removal of any $p-1$ vertices.
\hfill\break\noindent (2)
Let   $q>0$. Graph $G$ is
\emph{$(p,q)$-connected}, if $G$ remains connected after removal of any $p$ vertices and any $q-1$ edges. 
\end{definition}}

Clearly, if $G$ is $(p,q)$-connected graph, then $G$ is $(p^{\prime},q^{\prime})$-connected for any
$p^{\prime}\leq p$ and any $q^{\prime}\leq q$, 
and if $q>0$ then $(p,q)$-connected graph is also $(p+1,0)$-connected.
For any $(p,q)$-connected graph we have $p+q\leq\lambda(G)\leq\delta_G$, thus
each vertex of a $(p,q)$-connected graph has at least $p+q$ neighbors, 
and hence $(p,q)$-connected graph has at least $p+q+1$ vertices.

Mixed connectivity is  a generalization of vertex and  edge con\-nec\-ti\-vi\-ty:
a graph $G$ is  $(p,0)$-connected for all $p \leq \kappa(G)$ and is not $(p,0)$-connected for $p>  \kappa(G)$.
Furthermore, $G$ is  $(0,q)$-connected for all $q \leq \lambda(G)$ and is not $(0,q)$-connected for $q>  \lambda(G)$.
In particular, any graph $G$ is $(\kappa(G),0)$-connected and $(0,\lambda(G))$-connected.

%
%

The next statement follows directly from  Proposition \ref{T1}.

\begin{corollary} \label{T2}
If $p> 0$ and graph $G$ is $(p,q)$-connected then $G$ is $(p-1,q+1)$-connected.
\end{corollary}

Hence for $p>0$ we have a chain of implications 

\noindent
$(p,q)$-C $\Longrightarrow$ $(p-1,q+1)$-C $\Longrightarrow\dots\Longrightarrow$  $(1,p-1+q)$-C $\Longrightarrow$  $(0,p+q)$-C

\noindent
where $(i,j)$-C stands for "$G$ is $(i,j)$-connected".

Corollary \ref{T2} is for $q=0$  a generalization of well-known proposition 
that any $k$-connected graph is also $k$-edge connected. 
If $G$ is $(k,0)$-connected, then it is also $(k-i,i)$-connected for any $i\leq k$, and hence also $(0,k)$-connected.

If for graph $G$ $\kappa(G)=\lambda(G) =k$, then $G$ is $(i,j)$-connected exactly when $i+j\leq k$.
However, if  $2\leq \kappa(G)<\lambda(G)$, the question whether $G$ is $(i,j)$-connected  for  $1\leq i<\kappa(G)< i+j \leq \lambda(G)$
is not trivial. The example below shows that in general knowing 
$\kappa(G)$ and $\lambda(G)$ is not enough to decide whether  $G$ is  $(i,j)$-connected.

\begin{figure}[htb] 
\begin{center}
    \includegraphics[width=4.0in]{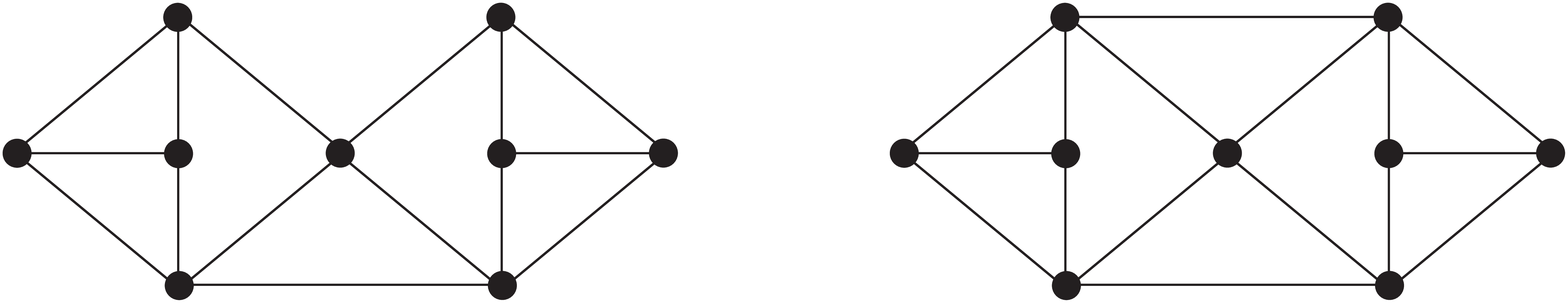}
    \caption{Graphs $G_1$ and $G_2$ from Example \ref{pr}.}
   \label{sl}
  \end{center}
\end{figure}

\begin{example} \label{pr}
For graphs on Fig. \ref{sl} we have $\kappa(G_1)=\kappa(G_2)=2$ and $\lambda(G_1)=\lambda(G_2)=3$. Both graphs are 
$(2,0)$-connected $\Longrightarrow$ $(1,1)$-connected $\Longrightarrow$ $(0,2)$-connected and $(0,3)$-connected.
Graph $G_1$ is not $(1,2)$-connected, while graph $G_2$ is.
\end{example}
 
Let $G$ be any connected graph, $ \kappa(G)<\lambda(G)$. All known mixed connectivities for graph $G$ are summarized 
in Diagram 1 (Fig. \ref{Diagram1}).

\begin{figure} 
\begin{scriptsize}
$$ \begin{array}{ccccccccccc}

  &   &   &   &      &      &      &       &       &            &  (0,\lambda(G)){\textrm -C}   \\  
  &   &   &   &      &      &      &       &       &            &  \Downarrow \\ 
  &   &   &   &      &      &      &       &       &						 &  (0,\lambda-1){\textrm -C} \\ 
  &   &   &   &      &      &      &       &       & 					 &  \Downarrow \\ 
  &   &   &	  &      &      &      &       &       & 					 &  \vdots \\
  &   &   &   &      &      &      &       &       & 					 &  \Downarrow \\ 
(\kappa(G),0){\textrm -C} & \Leftrightarrow    &  (\kappa-1,1){\textrm -C} & \Rightarrow  &  (\kappa-2,2){\textrm -C}    & \Rightarrow &  \dots  & \Rightarrow  & (1,\kappa-1){\textrm -C}    & \Rightarrow  &  (0,\kappa){\textrm -C} \\  
      &       & \Downarrow &   & \Downarrow &     &     &     & \Downarrow &      &  \Downarrow \\ 
      &       &  (\kappa-1,0){\textrm -C}  &  \Leftrightarrow &  (\kappa-2,1){\textrm -C} & \Rightarrow      & \dots & \Rightarrow  & (1,\kappa-2){\textrm -C}    & \Rightarrow  &  (0,\kappa-1){\textrm -C} \\  
      &       &  &   &\Downarrow &     & &            & \Downarrow &              &  \Downarrow \\ 
      &       &  &   &  (\kappa-2,0){\textrm -C}  &  \Leftrightarrow  \dots  \Rightarrow  & \dots & \Rightarrow  & (1,\kappa-3){\textrm -C}    & \Rightarrow  &  (0,\kappa-2){\textrm -C} \\  
   &   &  &        &  &  &  &    &  \Downarrow &       &  \Downarrow \\ 
   &   &  &       & & 							  &   &	 & \vdots &             & \vdots \\
   &   &  &      &  &  &  &   &  \Downarrow &       &  \Downarrow \\ 
   &   &  &      &  &  & (2,0){\textrm -C}  & \Leftrightarrow  & (1,1){\textrm -C}    & \Rightarrow  &  (0,2){\textrm -C} \\  
   &   &  &     & &   &       &       & \Downarrow &              &  \Downarrow \\ 
   &   &  &     & &   &       &       &  (1,0){\textrm -C}& \Leftrightarrow  &  (0,1){\textrm -C} \\ 

\end{array} 
$$
\end{scriptsize}
\caption{Diagram 1} \label{Diagram1}
\end{figure}

Both edge connectivity   and  connectivity (or vertex connectivity) of a  graph can be computed in 
polynomial time.  
More precisely, according to \cite{KorteVygen}, 
there is an algorithm with time complexity $O(n^4)$ for vertex connectivity \cite{Henzinger}
and edge connectivity can be computed within $O(m+\lambda(G) n \log\frac{n}{\lambda(G)})$ time \cite{Gabow}.
Hence given a graph $G$, there are polynomial time algorithms to decide $(i,0)$ and $(0,j)$-connectivities for all $i$, $j$.
Therefore it is interesting to ask 

\noindent{\bf Problem.}
Let $G$ be a graph and $1\leq i < \kappa(G) < i+j \leq \lambda(G)$.
Is there a polynomial algorithm to decide whether $G$ is $(i,j)$-connected?

%


If graph $G$ is $(i,j)$-connected for some $1\leq i<\kappa(G) < i+j \leq \lambda(G)$ then the upper part of 
diagram can be updated with known mixed connectivities as in Diagram 2 (Fig.\ref{Diagram2}).

\begin{figure} 
\begin{scriptsize}
$$ \begin{array}{ccccccccccccc}

  &   &   &   &   &&   &      &      &       &       &            &  (0,\lambda(G)){\textrm -C}   \\  
  &   &   &   &    &&  &      &      &       &       &            &  \Downarrow \\ 
  &   &   &   &    &&  &      &      &       &       &            &  \vdots \\ 
  &   &   &   &    &&  &      &      &       &       &            &  \Downarrow \\ 
   &   &  &   &&    & {\textbf (i,j){\textrm -C}}  &  \Rightarrow   & (i-1,j+1){\textrm -C}    &	\Rightarrow  & \dots  &   \Rightarrow           & (0,i+j){\textrm -C} \\
  &   &   &   &   && \Downarrow  &      &  \Downarrow     &       &       & 					 &  \Downarrow \\ 
   &   &  &   &&    & (i,j-1){\textrm -C} & 	 \Rightarrow  					  & (i-1,j){\textrm -C}    &	\Rightarrow  & \dots  &   \Rightarrow           & (0,i+j-1){\textrm -C} \\
  &   &   &   &   && \Downarrow  &      &  \Downarrow     &       &       & 					 &  \Downarrow \\ 
  &   &   &	  &   && \vdots  &      &  \vdots    &       &       & 					 &  \vdots \\
  &   &   &   &   && \Downarrow   &      &  \Downarrow     &       &       & 					 &  \Downarrow \\ 
(\kappa(G),0){\textrm -C} & \Leftrightarrow   & (\kappa-1,1){\textrm -C}   & \Rightarrow &\dots &\Rightarrow &  (i,\kappa-i){\textrm -C}    & \Rightarrow &  (i-1,\kappa-i+1){\textrm -C}  & \Rightarrow&  \dots    & \Rightarrow  &  (0,\kappa(G)){\textrm -C} \\  

\end{array} 
$$
\end{scriptsize}
\caption{Diagram 2} \label{Diagram2}  
~~\\
\end{figure}

Two extreme cases are:
\hfill\break\noindent
(1)
If graph $G$ is
$(\kappa(G)-1,\lambda(G) - \kappa(G)+1)$-connected then the diagram of mixed connectivities is maximal
because this connectivity implies all possible mixed connectivities for any connected graph. 
Namely, if graph $G$ is
$(\kappa(G)-1,\lambda(G) - \kappa(G)+1)$-connected then graph $G$ is $(i,j)$-connected 
for all $i<\kappa(G)$ and all $i+j\leq \lambda(G)$.
\hfill\break\noindent
(2)
If $G$ is not $(1,\kappa(G))$-connected then the diagram of mixed connectivities is minimal, i.e. the 
connectivities from Diagram 1 are all connectivities of $G$.

In general,   $(\kappa(G)-1)\times (\lambda(G) - \kappa(G))$ 
different connectivities  may have to be checked to complete the diagram of all connectivities of $G$.

\medskip

\section{Mixed connectivity of Cartesian graph bundles}

A Cartesian graph bundle is a generalization of graph cover and the Cartesian graph product.
Let $G_1$ and $G_2$ be graphs. The {\em Cartesian product} of graphs $G_1$ and $G_2$, $G=G_1\Box G_2$,
is defined on the vertex set $V(G_1)\times V(G_2)$.
Vertices $(u_1,v_1)$ and $(u_2,v_2)$ are 
adjacent if either $u_1u_2\in E(G_1)$ and $v_1=v_2$ or $v_1v_2\in E(G_2)$ and $u_1=u_2$. 
For further reading on graph products we recommend  \cite{3}. 

{\begin{definition}
Let $B$ and $F$ be graphs. A graph $G$ is a {\em Cartesian graph bundle with fibre $F$ over the base graph} $B$ if 
there is a {\em graph map} $p:G\rightarrow B$ such that for each vertex $v\in V(B)$, $p^{-1}(\{v\})$ is isomorphic to $F$, 
and for each edge $e=uv\in E(B)$, $p^{-1}(\{e\})$ is isomorphic to $F \Box K_2$.
\end{definition}} 
More precisely, the mapping $p:G\rightarrow B$  maps graph elements of $G$ to graph elements of $B$, 
i.e. $p: V(G) \cup E(G) \rightarrow  V(B) \cup E(B)$. In particular, here we also assume that 
the vertices of $G$ are mapped to vertices of $B$ and the edges of $G$ are mapped either to vertices or to edges of $B$. 
We say an edge $e\in E(G)$ is {\em degenerate} if $p(e)$ is a vertex. Otherwise we call it {\em nondegenerate}.
The mapping $p$ will also be called the {\em projection} (of the bundle $G$ to its base $B$).
Note that  each edge $e=uv \in E(B)$ naturally induces an isomorphism $\varphi_e :p^{-1}(\{u\})\rightarrow p^{-1}(\{v\})$
between two fibres. 
It may be interesting to note that while it is well-known that a graph can have only one representation  as a product
(up to isomorphism and up to the order of factors) \cite{3},
there may be many different graph bundle representations of the same graph \cite{ZmZe2002b}.
Here we assume that the bundle representation is given.
Note that in some cases finding a representation of $G$ as a graph bundle can be found in polynomial time 
\cite{ImPiZe,ZmZe2000,ZmZe2001a,ZmZe2002b,ZmZe2002a,directed}. 
For example, one of the easy classes are the Cartesian graph bundles over triangle-free base \cite{ImPiZe}.
Graph bundles were first studied in \cite{PiVr,ST_P_V}. 
Note that a graph bundle over a tree $T$ (as a base graph) with fibre $F$ is isomorphic to the 
Cartesian product $T\Box F$ (not difficult to see, appears already in  \cite{PiVr}), 
i.e. we can assume that all isomorphisms $\varphi_e$ are identities.
 
\begin{figure}[h!]
 \begin{center}
    \includegraphics[width=4in]{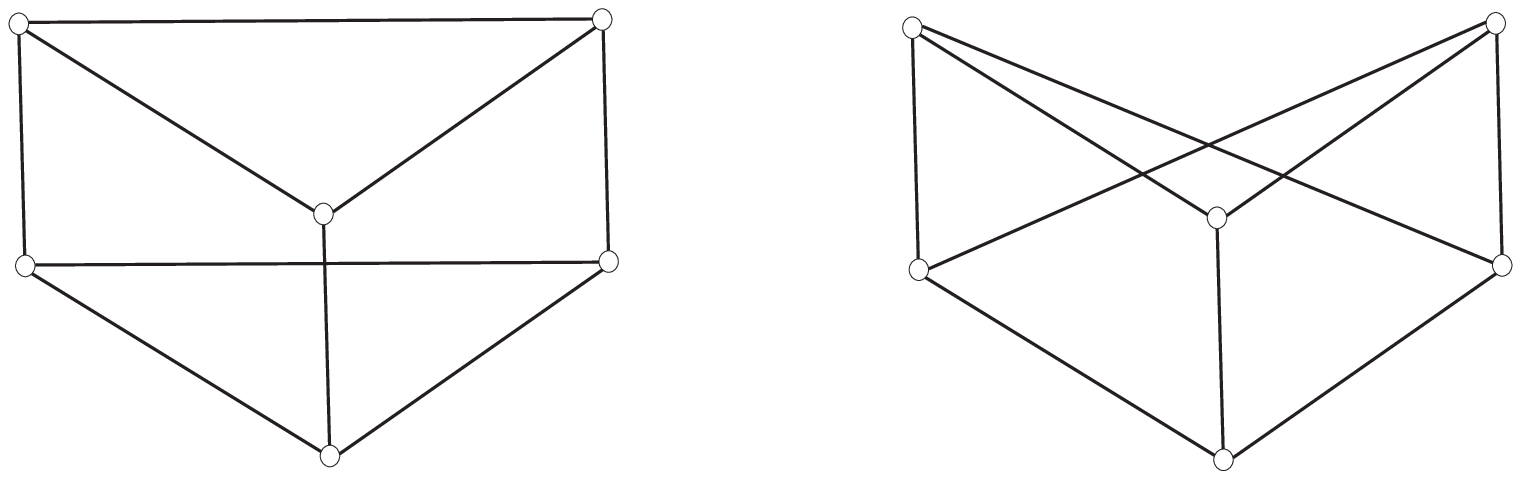}
    \caption{Nonisomorphic bundles from Example \ref{ey}}
    \label{bunddd}
  \end{center}
\end{figure}

\begin{example}
\label{ey}
Let $F=K_2$ and $B=C_3$. On Fig. \ref{bunddd} we see two nonisomorphic bundles with fibre $F$ over the base graph $B$.
Informally, one can say that bundles are "twisted products".
\end{example}

\begin{figure}[h!]
 \begin{center}
    \includegraphics[width=5.0in]{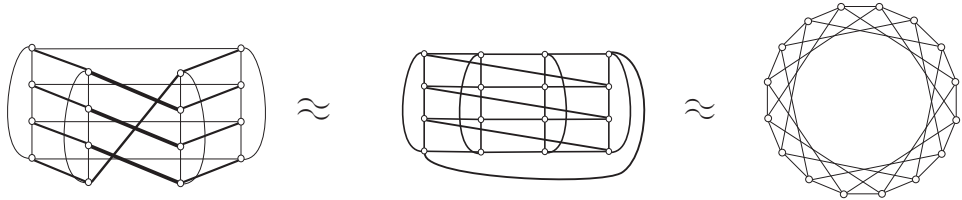}
    \caption{Twisted torus: Cartesian graph bundle with fibre $C_4$ over base $C_4$.}
    \label{iliac}
  \end{center}
\end{figure}

\begin{example}
It is less known that graph bundles also appear as computer topologies. 
A well known example is the twisted torus on Fig. \ref{iliac}.
Cartesian graph bundle with fibre $C_4$ over base $C_4$ is the ILLIAC IV architecture \cite{ILLIAC},
a famous supercomputer that inspired some modern multicomputer architectures.
It may be interesting to note that the original design was a graph bundle with fibre $C_8$ over base $C_8$, 
but due to high cost a smaller version was build \cite{computermuseum}.
\end{example}


Let $G$ be a Cartesian graph bundle with fibre $F$ over the base graph $B$.
The {\em fibre of vertex} $x \in V(G)$ is denoted by $F_x$, formally, $F_x = p^{-1}(\{p(x)\})$. We will also use notation $F(u)$ 
for the fibre of the vertex $u \in V(B)$, i.e. $F(u) = p^{-1}(\{u\})$. Note that $F_x = F(p(x))$. 
Let $u,v\in V(B)$ be distinct vertices, 
$Q$ be a path from $u$ to $v$ in $B$, and $x\in F(u)$. Then the {\em lift of the path $Q$ to the vertex} $x\in V(G)$, 
$\tilde Q_x$, is the path from $x\in F(u)$ to a vertex in $F(v)$, 
such that $p(\tilde Q_x)=Q$ and $\ell(\tilde Q_x)=\ell(Q)$.
Let $x,x'\in F(u)$. Then $\tilde Q_x$ and $\tilde Q_{x'}$ have different endpoints in $F(v)$ and are disjoint paths if 
and only if $x\neq x'$.
We will also use notation $\tilde Q$ for lifts of path $Q$ to any vertex in $F(u)$. 

In previous  work \cite{zer-ban,erves} 
on vertex and edge fault diameters of Cartesian graph bundles  propositions 
about vertex and edge con\-nec\-ti\-vi\-ty of Cartesian graph bundles have been proved. 
In terms of mixed connectivity they read as follows. 

\begin{proposition} \cite{zer-ban}
\label{pt} If graph $F$ is $(p_{F},0)$-connected and graph $B$ is $(p_{B},0)$-connected, then 
Cartesian graph bundle with fibre $F$ over the base graph $B$ is $(p_{F}+p_{B},0)$-connected.
\end{proposition}

\begin{proposition} \cite{erves}
\label{pp} If graph $F$ is $(0,q_{F})$-connected and graph $B$ is $(0,q_{B})$-con\-nected, then 
Cartesian graph bundle with fibre $F$ over the base graph $B$ is $(0,q_{F}+q_{B})$-connected.
\end{proposition}

A natural generalization of these two propositions would be that 
if graph $F$ is $(p_{F},q_{F})$-connected and graph $B$ is $(p_{B},q_{B})$-connected, then 
Cartesian graph bundle with fibre $F$ over the base graph $B$ is $(p_{F}+p_{B},q_{F}+q_{B})$-connected.
This is indeed true, but we can prove a slightly stronger statement (see Theorem \ref{mpsv}).
Roughly speaking, 
consider the maximum allowed number of faulty elements of graphs $F$, $B$ and $G$. 
If $q_{F}=0$ or $q_{B}=0$, then generalization assures connectivity of graph bundle $G$ with one more faulty element 
(vertex or edge) as the sum of faulty vertices and faulty edges in graphs $F$ and $B$. 

We prove that besides the sum of allowed  faulty elements of the fibre and the base, 
one additional faulty element  is allowed.
We also show that whenever applicable, the extra faulty element can be a node.
In particular this improves the result on edge connectivity.  
Namely, Theorem \ref{mpsv} is for $p_F=p_B=0$ stronger than Proposition \ref{pp}. 
If graph $F$ is $(0,q_{F})$-connected and graph $B$ is $(0,q_{B})$-con\-nected, then 
by Theorem \ref{mpsv}, Cartesian graph bundle with fibre $F$ over the base graph $B$ is 
$(1,q_{F}+q_{B}-1)$-connected,
while Proposition \ref{pp} assures only 
$(0,q_{F}+q_{B})$-connectivity.


\begin{theorem}
\label{mpsv} Let $G$ be a Cartesian graph bundle with fibre $F$ over the base graph $B$, graph $F$ be $(p_{F},q_{F})$-connected and
graph $B$ be $(p_{B},q_{B})$-connected.
Then Cartesian graph bundle $G$ is:

\begin{enumerate}
\item $(p_{F}+p_{B},q_{F}+q_{B})$-connected.
\item for $q_{F},q_{B}>0$ also $(p_{F}+p_{B}+1,q_{F}+q_{B}-1)$-connected.
\end{enumerate}
\end{theorem}

As the Cartesian product is a Cartesian graph bundle where all the isomorphisms between the fibres are identities, 
the statement about mixed connectivity of Cartesian graph products of a finite number of factors follows 
easily from Theorem \ref{mpsv}. 
Let 
$G=G_1\Box G_2\Box \dots \Box G_k$, and let $G_i, i=1,\dots,k,$ be $(p_i,q_i)$-connected, and $q_i>0$. 
Then, by induction, $G$ is $(\sum p_i+k-1,\sum q_i-k+1)$-connected.
Therefore

\begin{corollary}
\label{mppr} Let graphs $G_i, i=1,\dots,n,$ be $(p_i,q_i)$-connected and let $k$ be the number of graphs with $q_i>0$. 
Then the Cartesian graph product $G=G_1\Box G_2\Box \dots \Box G_n$ is:

\begin{enumerate}
\item $(\sum p_i,\sum q_i)$-connected, and
\item $(\sum p_i+k-1,\sum q_i-k+1)$-connected, for $ k\geq 1$. 
\end{enumerate}
\end{corollary}

\section{Proof of the main theorem}

Let $G$ be any connected graph. Then for any $p<\lambda(G)$ a graph $G$ is 
$(p+1,0)$-connected if and only if $G$ is $(p,1)$-connected.
Therefore it is enough to prove Theorem \ref{mpsv} only for the case when $q_{F},q_{B}>0$.
For example, if graphs $F$ and $B$ are $(p_{F},0)$-connected and
$(p_{B},0)$-connected, then graphs are also $(p_{F}-1,1)$-connected and $(p_{B}-1,1)$-connected respectively.
By Theorem \ref{mpsv} $\textit{(2)}$ a Cartesian graph bundle $G$ is $(p_{F}+p_{B}-1,1)$-connected, 
hence $G$ is $(p_{F}+p_{B},0)$-connected.

Hence from now on we can assume $q_{F},q_{B}>0$.
We will show that Cartesian graph bundle with fibre $F$ over the base graph $B$ without maximum allowed number of  
faulty elements is  connected.
Denote the set of faulty vertices by  $X\subseteq V(G)$, $\left\vert X\right\vert =p_{F}+p_{B}+1$, and the set of 
faulty edges by  $Y\subseteq E(G)$, $\left\vert Y\right\vert =q_{F}+q_{B}-2$.
We have to show that $G\setminus(X\cup Y)$ is connected graph.

For each vertex $v\in V(B)$, fibre $F(v)\setminus(X\cup Y)$ is either connected or disconnected subgraph of $G\setminus(X\cup Y)$.
This wording is convenient and we will use it although 
it is not formally clean  because $G\setminus(X\cup Y)$ is of course most likely not a graph bundle, 
and  by saying  that the fibre $F(v)$ is connected or disconnected  we are in fact referring to 
the properties of the subgraph $F(v)\setminus(X\cup Y)$.

It is not difficult to see that a graph $G\setminus(X\cup Y)$ always contains at least one connected fibre, 
and it may contain some disconnected fibres.  
We will prove Theorem \ref{mpsv} by proving  two lemmas.   
With Lemma \ref{lema1} we show that in $G\setminus(X\cup Y)$ any vertex of a disconnected fibre is connected 
with some vertex of connected fibre.
With Lemma \ref{lema2} we show that there is a path between any two connected fibres in $G\setminus(X\cup Y)$.
Both lemmas together assure that all pairs of vertices in $G\setminus(X\cup Y)$ are connected by paths, which 
implies Theorem \ref{mpsv}.

 The {\em weights} which are used in proofs of lemmas are defined as follows:
\begin{itemize}
\item {\em (faulty) vertex-weight of vertex} $v\in V(B)$, $w_X(v)$, 
      is the  number of faulty vertices in fibre $F(v)$, $w_X(v)=\left\vert F(v)\cap X\right\vert$; 
\item {\em (faulty) edge-weight of vertex} $v\in V(B)$, $w_Y(v)$, 
      is the number of faulty (degenerate) edges in fibre $F(v)$, $w_Y(v)=\left\vert F(v)\cap Y\right\vert$; 
\item {\em (faulty) edge-weight of path} $Q\subseteq B$, $w_Y(Q)$, 
      is the number of faulty (nondegenerate) edges on lifts of path $Q$, $w_Y(Q)=\left\vert p^{-1}(Q)\cap Y_N\right\vert$, 
      where $Y_N$  is the set of faulty nondegenerate edges in $Y$, $Y_N\subseteq Y$, $p(Y_N)\subseteq E(B)$, $p(Y\setminus Y_N)\subseteq V(B)$.
\end{itemize}

\begin{lemma}
\label{lema1} Let $G$ be a Cartesian graph bundle with fibre $F$ over the base graph $B$, 
graph $F$ be $(p_{F},q_{F})$-connected, graph $B$ be  $(p_{B},q_{B})$-connected.
Let $X\subseteq V(G)$ and $Y\subseteq E(G)$ be sets of faulty vertices and faulty edges
with maximal allowed number of elements, and let $x\in V(G)\setminus X$ be any vertex of a disconnected fibre in $G\setminus(X\cup Y)$. 
Then in $G\setminus(X\cup Y)$ exists (neighboring) connected fibre and there is a path between vertex $x$ and a vertex of connected fibre. 
\end{lemma}

\proof{Proof.}
Let graph $F$ be $(p_{F},q_{F})$-connected, graph $B$ be $(p_{B},q_{B})$-connected, $q_{F},q_{B}>0$,
let $X\subseteq V(G)$ and $Y\subseteq E(G)$ be sets of faulty vertices and faulty edges, $\left\vert X\right\vert =p_{F}+p_{B}+1$, $\left\vert Y\right\vert=q_{F}+q_{B}-2$. 
Let $x\in V(G)\setminus X$ and $F_x\setminus(X\cup Y)$ be disconnected fibre in $G\setminus(X\cup Y)$.
Then either $w_X(p(x))> p_{F}$ or $w_X(p(x))+w_Y(p(x))\geq p_{F}+q_{F}$. 
We distinguish two cases.

\begin{enumerate}
\item Suppose $w_X(p(x))+w_Y(p(x))\geq p_{F}+q_{F}$.

\noindent Fibre $F_x$ contains at least $p_F+q_F$ faulty elements, so outside of fibre $F_x$ there are at 
most $p_B+q_B-1$ faulty elements. 
As vertex $p(x)$ has at least $p_B+q_B$ neighbors in $B$, there is a neighbor $v$, $e=p(x)v$, 
with weights $w_X(v)+w_Y(v)=0$ and $w_Y(e)=0$. Hence fibre $F(v)\subseteq G\setminus(X\cup Y)$ is 
connected  and the lift $\tilde e_x$ avoids $X\cup Y$. So $x$ is adjacent to a vertex in $F(v)$ as needed. 

\item Now assume that  $w_X(p(x))> p_{F}$. 

\noindent Fibre $F_x$ contains at least $p_F+1$ faulty vertices, so outside of fibre $F_x$ there are at 
most $p_B$ faulty vertices (and $q_{F}+q_{B}$ faulty edges). 

\noindent Let $X_{B}=\{  v\in V(B)\setminus\{  p(x)\}
;w_X(v) >0\}  $ and $b=\left\vert
X_{B}\right\vert $. 
Then $b\leq p_B$.

\noindent As there are at least $p_B+q_B$ neighbors of vertex $p(x)$ in $B$, 
therefore there are at least $p_B+q_B-b\geq 1$ neighbors $v_i$ in $B\setminus X_B$ with weights $w_X(v_i)=0$. 

\noindent If there is a neighbor $v$ of vertex $p(x)$ in $B\setminus X_B$ ($v=v_i$ for some $i$), 
$e=p(x)v$, with edge-weights $w_Y(v)=w_Y(e)=0$, then $F(v)\subseteq G\setminus(X\cup Y)$ is connected 
fibre and lift $\tilde e_x$ avoids $X\cup Y$. So $x$ is adjacent to a vertex in $F(v)$ as needed.

\noindent 
Now suppose that for every neighbor $v_i$ of vertex $p(x)$ in $B\setminus X_B$, $e_i=p(x)v_i$, $w_Y(v_i)+w_Y(e_i)>0$. 
Then let $v$ be any neighbor ($v=v_i$ for any $i$), $e=p(x)v$, in $B\setminus X_B$. Outside of $p^{-1}(\{e\})$ there
are at least $b$ faulty vertices (because we eliminate set $X_B$) and $p_B+q_B-b-1$ faulty edges (in other neighbors), 
together $p_B+q_B-1$ faulty  elements.
Therefore there are at most $p_{F}+q_{F}$ faulty elements in $p^{-1}(\{e\})$. As fibre $F_x\subset p^{-1}(\{e\})$ 
has at least $p_F+1$ faulty vertices, there are at most $q_{F}-1$ faulty edges in $F(v)\subset p^{-1}(\{e\})$, 
hence $F(v)\setminus Y\subseteq G\setminus(X\cup Y)$ is connected fibre. 

\noindent There are at least $p_F+q_F$ neighbors of vertex $x$ in $F_x$. Denote the neighbors 
by  $s_{i}$, $i=1,\dots,p_F+q_F$, and let $e_{i}=xs_{i}$.
There are $p_F+q_F+1$ vertex disjoint paths in $p^{-1}(\{e\})$ with one endpoint $x$ and another endpoint in 
connected fibre $F(v)$: $x,e_i,s_i,\tilde e,v_i$ and $x,\tilde e,v^{\prime}$, where $\tilde e$ is a lift of the edge 
$e=p(x)v$ and $v_{i},v^{\prime}$ are different vertices in fibre $F(v)$. As there are more vertex disjoint paths than 
faulty elements, at least one of these paths avoids $X\cup Y$. \qed
\end{enumerate}

\begin{lemma}
\label{lema2} Let $G$ be a Cartesian graph bundle with fibre $F$ over the base graph $B$, 
graph $F$ be $(p_{F},q_{F})$-connected, graph $B$ be $(p_{B},q_{B})$-connected. Let
$X\subseteq V(G)$ and $Y\subseteq E(G)$ be sets of faulty vertices and faulty edges with maximal allowed number of elements. 
Then for any two connected fibres $F_x\setminus(X\cup Y)$ and $F_y\setminus(X\cup Y)$ there is a path with 
endpoints in fibres $F_x$ and $F_y$ that avoids faulty elements.
\end{lemma}

\proof{Proof.}
Let graph $F$ be $(p_{F},q_{F})$-connected, graph $B$ be $(p_{B},q_{B})$-connected, $q_{F},q_{B}>0$,
let $X\subseteq V(G)$ and $Y\subseteq E(G)$ be sets of faulty vertices and faulty edges, $\left\vert X\right\vert =p_{F}+p_{B}+1$, $\left\vert Y\right\vert=q_{F}+q_{B}-2$. 
\newline
Let $X_{B}=\{  v\in V(B)\setminus\{  p(x),p(y)\}
;w_X(v) >0\}  $ and $b=\left\vert
X_{B}\right\vert $. 

\begin{enumerate}
\item First assume $b< p_B$.

\noindent Graph $B\setminus X_B$ is $(p_B-b,q_{B})$-connected, hence $B\setminus X_B$ is $(0,p_B+q_{B}-b)$-connected, 
and therefore there are at least $p_B+q_{B}-b\geq q_{B}+1\geq 2$ edge disjoint paths between 
$p(x)$ and $p(y)$ in $B$ that avoid $X_B$; one of them may be the edge. 
There are at least $p_B+q_{B}-b-1\geq q_{B}\geq1$ edge disjoint paths (with lengths more than $1$) that internally avoid $p(X)$.

\begin{itemize}
\item[(A)] Suppose there is a path $Q$ between $p(x)$ and $p(y)$ (with length more than $1$) that 
internally avoids $p(X)$ with edge-weights $w_Y(Q)+w_Y(v)=0$, where $v\in Q$ is neighbor of vertex $p(x)$. 
Then $F(v) \subseteq G\setminus(X\cup Y)$ is connected fibre. There are at least $p_F+q_F+1$ lifts  
of edge $e=p(x)v\subset Q$ with different endpoints in fibres  $F_x$ and $F(v)$. As these lifts 
contain at most $p_{F}$ faulty endpoints in fibre $F_x$,
there is a lift $\tilde e$ that avoids faulty elements. 
Similarly, there is a lift of the path $vp(y)\subset Q$ that avoids faulty elements, and as   
$F(v)\subseteq G\setminus(X\cup Y)$ is connected fibre, there is a path between fibres $F_x$ and $F_y$ that avoids faulty elements.

\item[(B)] 
Now suppose that for each path $Q_i$ between $p(x)$ and $p(y)$ (with length more than $1$) 
that internally avoids $p(X)$, the sum of edge-weights is $w_Y(Q_i)+w_Y(v_i)>0$, 
where $v_i\in Q_i$ is the neighbor of vertex $p(x)$ along $Q_i$. 
Let $Q$ be one of the shortest paths between $p(x)$ and $p(y)$ in $B\setminus X_B$. Then outside $p^{-1}(Q)$ 
there are at least $b$ faulty vertices and   $p_B+q_{B}-b-1$ faulty edges, 
together at least $p_B+q_{B}-1$ faulty elements. 
Therefore $p^{-1}(Q)$ contains at most $p_F+q_F$ faulty elements. As there are at least $p_F+q_F+1$ lifts of path $Q$ 
between fibres $F_x$ and $F_y$, at least one of them avoids faulty elements. 
\end{itemize}

\item Now assume $b\geq p_B$.

\noindent Let $X_B^{\prime} \subseteq X_{B}$ be any subset with $\left\vert X_{B}^{\prime}\right\vert =p_{B}$. Graph $B\setminus X_B^{\prime}$ is $(0,q_{B})$-connected, so there are $q_{B}\geq 1$ edge disjoint paths between $p(x)$ and $p(y)$ in $B\setminus X_B^{\prime}$.

\begin{itemize}
\item[(A)] If there is a path $Q$ between $p(x)$ and $p(y)$ in $B\setminus X_B^{\prime}$ with edge-weight $w_Y(Q)=0$, 
then the lifts of path $Q$ do not contain faulty edges. 
Outside $p^{-1}(Q)$ there are at least $p_B$ faulty vertices. Therefore $p^{-1}(Q)$ contains at most $p_F+1$ faulty vertices. As there are at least $p_F+q_F+1$ 
lifts of path $Q$ between fibres $F_x$ and $F_y$, at least one of them avoids faulty elements.

\item[(B)] Finally assume that there is no path with edge-weight $0$ between $p(x)$ and $p(y)$ in $B\setminus X_B^{\prime}$. 
Let $Q$ be any path between $p(x)$ and $p(y)$ in $B\setminus X_B^{\prime}$. Outside $p^{-1}(Q)$ there are at least 
$p_B$ faulty vertices and $q_{B}-1$ faulty edges, together at least $p_B+q_{B}-1$ faulty elements. 
Therefore $p^{-1}(Q)$ contains at most  $p_F+q_F$ faulty elements. 
As there are at least $p_F+q_F+1$ lifts of path $Q$ between fibres $F_x$ and $F_y$, 
at least one of them avoids faulty elements. \qed
\end{itemize}
\end{enumerate}


\end{document}